\documentclass[a4paper,12pt]{article}

\usepackage{amsthm}
\usepackage{amsmath,amssymb,latexsym,amsfonts}

\theoremstyle{plain}
\newtheorem{Thm}{Theorem}[section]
\newtheorem{Lem}[Thm]{Lemma}

\theoremstyle{definition}
\newtheorem{Def}[Thm]{Definition}

\newcommand{\Proof}[2][Proof]{\begin{proof}[{#1}] #2 \end{proof}}

\renewcommand{\hat}{\widehat}

\newcommand{\eps}{\ensuremath{\varepsilon}}

\newcommand{\dist}{\stackrel{{\rm d}}{=}}
\newcommand{\tend}[2]{\mathrel{\mathop{\longrightarrow}\limits^{#1}_{#2}}}

\newcommand{\rbra}[1]{\left( #1 \right)} 
\newcommand{\cbra}[1]{\left\{ #1 \right\}} 
\newcommand{\abra}[1]{\left\langle #1 \right\rangle} 

\newcommand{\bmat}[1]{\begin{bmatrix} #1 \end{bmatrix}} 

\newcommand{\bN}{\ensuremath{\mathbb{N}}}
\newcommand{\bR}{\ensuremath{\mathbb{R}}}
\newcommand{\bZ}{\ensuremath{\mathbb{Z}}}

\numberwithin{equation}{section}

\newcommand{\Trans}{Q}
\newcommand{\trans}{q}

\newcommand{\minip}[2]{
\begin{minipage}{#1}
\begin{center}
#2
\end{center}
\end{minipage}}

\begin{document}
\begin{center}
{\Large \bf Realization of finite-state mixing Markov chain 
as a random walk subject to a synchronizing road coloring 
}
\end{center}

\begin{center}
Kouji Yano\footnote{
Graduate School of Science, Kobe University, Kobe, JAPAN.}\footnote{
The research of this author was supported by KAKENHI (20740060).} 
\quad and \quad 
Kenji Yasutomi\footnote{
Department of Mathematical Sciences, 
Ritsumeikan University, Kusatsu, JAPAN.} 
\end{center}

\begin{center}
{\small \today}
\end{center}

\begin{abstract}
A mixing Markov chain is proved to be realized as a random walk in a directed graph 
subject to a synchronizing road coloring. 
The result ensures existence of appropriate random mappings 
in Propp--Wilson's coupling from the past. 
The proof is based on the road coloring theorem. 
A necessary and sufficient condition for approximate preservation of entropies 
is also given. 
\end{abstract}

\noindent
{\footnotesize Keywords and phrases: Markov chain, 
random walk in a directed graph, road coloring problem, 
Tsirelson's equation, coupling from the past.} 
\\
{\footnotesize AMS 2010 subject classifications: 
Primary
60J10; 
secondary
37A35; 
05C81; 
37H10. 
}

\section{Introduction}

Our purpose is to realize a mixing Markov chain 
as a suitable random walk in a directed graph, 
which is generated by a sequence of 
independent, identically distributed ({\em IID} for short) random variables 
taking values in the set of mappings of the state space.

\subsection{Notations}

Let $ V $ be a set of finite symbols, say $ V = \{ 1,\ldots,m \} $. 
Let $ Y=(Y_k)_{k \in \bZ} $ be a Markov chain taking values in $ V $ 
and indexed by $ \bZ $, the set of all integers. 
We call $ Y $ {\em stationary} if $ (Y_{k+n})_{k \in \bZ} \dist (Y_k)_{k \in \bZ} $ 
for all $ n \in \bZ $. 
We write $ \Trans = (\trans_{x,y})_{x,y \in V} $ 
for the one-step transition probability matrix of $ Y $, 
i.e., 
\begin{align}
\trans_{x,y} = P(Y_1=y|Y_0=x) 
, \quad x,y \in V . 
\label{}
\end{align}
The $ n $-th transition probability matrix is given 
by the $ n $-th product $ \Trans^n = (\trans^n_{x,y})_{x,y \in V} $. 
We call $ Y $ {\em irreducible} if 
for any $ x,y \in V $ there exists a positive number $ n=n(x,y) $ 
such that $ \trans^n_{x,y} >0 $. 
We call $ Y $ {\em aperiodic} if 
the greatest common divisor among $ \{ n \ge 1 : \trans^n_{x,x} > 0 \} $ is one 
for all $ x \in V $. 
We call $ Y $ {\em mixing} 
if $ Y $ is both irreducible and aperiodic, 
which is equivalent to the condition that 
there exists a positive integer $ r $ 
such that $ \trans^r_{x,y} > 0 $ for all $ x,y \in V $. 

Let $ \Sigma $ denote the set of all mappings from $ V $ to itself. 
For $ \sigma_1,\sigma_2 \in \Sigma $ and $ x \in V $, 
we write $ \sigma_2 \sigma_1 x $ simply for $ \sigma_2(\sigma_1(x)) $. 

\begin{Def}
Let $ \Trans $ be the one-step transition probability matrix 
of a stationary Markov chain. 
A probability law $ \mu $ on $ \Sigma $ is called a {\em mapping law} for $ \Trans $ if 
\begin{align}
\trans_{x,y} = \sum_{\sigma \in \Sigma : \sigma x = y} \mu(\sigma) 
, \quad x,y \in V . 
\label{eq: p1 and mu}
\end{align}
\end{Def}

\begin{Def}
For a probability law $ \mu $ on $ \Sigma $, 
a {\em $ \mu $-random walk} is a Markov chain 
$ (X,N) = (X_k,N_k)_{k \in \bZ} $ taking values in $ V \times \Sigma $ such that 
$ N=(N_k)_{k \in \bZ} $ is IID with common law $ \mu $ 
such that each $ N_k $ is independent of $ \sigma(X_j,N_j:j \le k-1) $ 
and 
\begin{align}
X_k = N_k X_{k-1} 
\quad \text{a.s. for $ k \in \bZ $}. 
\label{eq: Tsirel}
\end{align}
\end{Def}

Let $ Y=(Y_k)_{k \in \bZ} $ be a stationary Markov chain 
with one-step transition probability matrix $ \Trans $. 
Let $ (X,N) $ be a $ \mu $-random walk. 
Then it is obvious that 
$ Y \dist X $ if and only if 
$ \mu $ is a mapping law for $ \Trans $. 
For any stationary Markov chain $ Y $, 
we can find a mapping law $ \mu $ for $ \Trans $ 
(see Lemma \ref{lem: rational}).

Let us illustrate our notations. 
See Figure 1 below, where $ V=\{ 1,2,3 \} $ and 
\begin{align}
\Trans = \bmat{ 
q_{1,1} & q_{1,2} & q_{1,3} \\
q_{2,1} & q_{2,2} & q_{2,3} \\
q_{3,1} & q_{3,2} & q_{3,3} }
= \bmat{ 
  0 & 2/3 & 1/3 \\
1/3 &   0 & 2/3 \\
2/3 & 1/3 &   0 }. 
\label{}
\end{align}
Let $ \sigma^{(1)} $, $ \sigma^{(2)} $, $ \sigma^{(3)} $ and $ \sigma^{(4)} $ 
be elements of $ \Sigma $ characterized by Figures 2, 3, 4 and 5 below, respectively. 
The transition probability $ \Trans $ possesses several mapping laws; 
among others, we have $ \mu^{(1)} $ and $ \mu^{(2)} $ defined as follows: 
\begin{align}
\mu^{(1)}(\sigma^{(1)}) = \mu^{(1)}(\sigma^{(2)}) = \mu^{(1)}(\sigma^{(3)}) = 1/3 
, \label{eq: mu1} \\
\mu^{(2)}(\sigma^{(3)}) = 2/3 , \quad \mu^{(2)}(\sigma^{(4)}) =  1/3 
. \label{eq: mu2}
\end{align}
The two random walks $ (X,N) $ corresponding to $ \mu^{(1)} $ and $ \mu^{(2)} $ 
have distinct joint laws, 
but have identical marginal law of $ X $ 
which is a Markov chain with one-step transition probability $ \Trans $. 

\begin{center}
\unitlength 0.1in
\begin{picture}( 20.0000, 16.0000)(  0.0000,-16.0000)
%
\special{pn 8}%
\special{ar 1000 200 200 200  0.0000000 6.2831853}%
%
\special{pn 8}%
\special{ar 200 1400 200 200  0.0000000 6.2831853}%
%
\special{pn 8}%
\special{ar 1800 1400 200 200  0.0000000 6.2831853}%
%
\special{pn 8}%
\special{pa 830 330}%
\special{pa 250 1200}%
\special{fp}%
\special{pa 250 1060}%
\special{pa 250 1210}%
\special{fp}%
\special{pa 250 1210}%
\special{pa 380 1150}%
\special{fp}%
\special{pa 390 1480}%
\special{pa 1610 1480}%
\special{fp}%
\special{pa 1480 1400}%
\special{pa 1600 1460}%
\special{fp}%
\special{pa 1600 1460}%
\special{pa 1470 1580}%
\special{fp}%
\special{pa 1630 1260}%
\special{pa 1040 400}%
\special{fp}%
\special{pa 1050 550}%
\special{pa 1050 400}%
\special{fp}%
\special{pa 1050 400}%
\special{pa 1200 450}%
\special{fp}%
%
\special{pa 380 1300}%
\special{pa 950 400}%
\special{fp}%
\special{pa 800 470}%
\special{pa 940 400}%
\special{fp}%
\special{pa 940 400}%
\special{pa 940 580}%
\special{fp}%
\special{pa 1600 1320}%
\special{pa 400 1320}%
\special{fp}%
\special{pa 520 1230}%
\special{pa 390 1310}%
\special{fp}%
\special{pa 390 1310}%
\special{pa 530 1390}%
\special{fp}%
\special{pa 1170 320}%
\special{pa 1780 1200}%
\special{fp}%
\special{pa 1650 1150}%
\special{pa 1790 1200}%
\special{fp}%
\special{pa 1790 1200}%
\special{pa 1800 1070}%
\special{fp}%
\put(10.0000,-2.2000){\makebox(0,0){1}}%
\put(1.9000,-14.1000){\makebox(0,0){2}}%
\put(18.0000,-14.0000){\makebox(0,0){3}}%
\put(4.6000,-6.5000){\makebox(0,0){$ 2/3 $}}%
\put(15.6000,-6.5000){\makebox(0,0){$ 1/3 $}}%
\put(7.5000,-9.6000){\makebox(0,0){$ 1/3 $}}%
\put(12.5000,-9.6000){\makebox(0,0){$ 2/3 $}}%
\put(10.0000,-12.2000){\makebox(0,0){$ 1/3 $}}%
\put(10.0000,-16.0000){\makebox(0,0){$ 2/3 $}}%
\end{picture}%
\\[2mm] Figure 1. Transition probability
\end{center}

\begin{center}
\minip{7cm}{
\unitlength 0.1in
\begin{picture}( 20.0000, 16.0000)(  0.0000,-16.0000)
%
\special{pn 8}%
\special{ar 1000 200 200 200  0.0000000 6.2831853}%
%
\special{pn 8}%
\special{ar 200 1400 200 200  0.0000000 6.2831853}%
%
\special{pn 8}%
\special{ar 1800 1400 200 200  0.0000000 6.2831853}%
%
\special{pn 8}%
\special{pa 390 1480}%
\special{pa 1610 1480}%
\special{fp}%
\special{pa 1480 1400}%
\special{pa 1600 1460}%
\special{fp}%
\special{pa 1600 1460}%
\special{pa 1470 1580}%
\special{fp}%
\special{pa 1630 1260}%
\special{pa 1040 400}%
\special{fp}%
\special{pa 1050 550}%
\special{pa 1050 400}%
\special{fp}%
\special{pa 1050 400}%
\special{pa 1200 450}%
\special{fp}%
\special{pa 1170 320}%
\special{pa 1780 1200}%
\special{fp}%
\special{pa 1650 1150}%
\special{pa 1790 1200}%
\special{fp}%
\special{pa 1790 1200}%
\special{pa 1800 1070}%
\special{fp}%
\put(10.0000,-2.2000){\makebox(0,0){1}}%
\put(1.9000,-14.1000){\makebox(0,0){2}}%
\put(18.0000,-14.0000){\makebox(0,0){3}}%
\end{picture}%
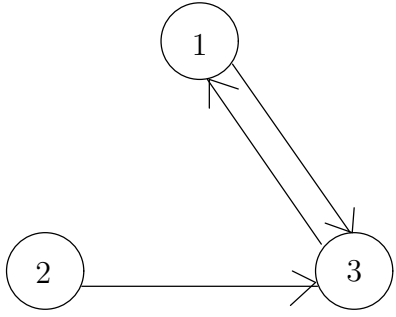
\\ Figure 2. $ \sigma^{(1)} $
}
\minip{7cm}{
\unitlength 0.1in
\begin{picture}( 20.0000, 16.0000)(  0.0000,-16.0000)
%
\special{pn 8}%
\special{ar 1000 200 200 200  0.0000000 6.2831853}%
%
\special{pn 8}%
\special{ar 200 1400 200 200  0.0000000 6.2831853}%
%
\special{pn 8}%
\special{ar 1800 1400 200 200  0.0000000 6.2831853}%
%
\special{pn 8}%
\special{pa 830 330}%
\special{pa 250 1200}%
\special{fp}%
\special{pa 250 1060}%
\special{pa 250 1210}%
\special{fp}%
\special{pa 250 1210}%
\special{pa 380 1150}%
\special{fp}%
\special{pa 380 1300}%
\special{pa 950 400}%
\special{fp}%
\special{pa 800 470}%
\special{pa 940 400}%
\special{fp}%
\special{pa 940 400}%
\special{pa 940 580}%
\special{fp}%
\special{pa 1600 1320}%
\special{pa 400 1320}%
\special{fp}%
\special{pa 520 1230}%
\special{pa 390 1310}%
\special{fp}%
\special{pa 390 1310}%
\special{pa 530 1390}%
\special{fp}%
\put(10.0000,-2.2000){\makebox(0,0){1}}%
\put(1.9000,-14.1000){\makebox(0,0){2}}%
\put(18.0000,-14.0000){\makebox(0,0){3}}%
\end{picture}%
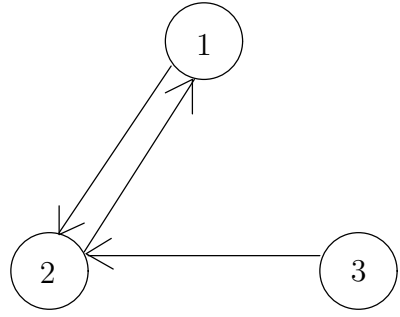
\\ Figure 3. $ \sigma^{(2)} $
}
\end{center}
\begin{center}
\minip{7cm}{
\unitlength 0.1in
\begin{picture}( 20.0000, 16.0000)(  0.0000,-16.0000)
%
\special{pn 8}%
\special{ar 1000 200 200 200  0.0000000 6.2831853}%
%
\special{pn 8}%
\special{ar 200 1400 200 200  0.0000000 6.2831853}%
%
\special{pn 8}%
\special{ar 1800 1400 200 200  0.0000000 6.2831853}%
%
\special{pn 8}%
\special{pa 830 330}%
\special{pa 250 1200}%
\special{fp}%
\special{pa 250 1060}%
\special{pa 250 1210}%
\special{fp}%
\special{pa 250 1210}%
\special{pa 380 1150}%
\special{fp}%
\special{pa 390 1480}%
\special{pa 1610 1480}%
\special{fp}%
\special{pa 1480 1400}%
\special{pa 1600 1460}%
\special{fp}%
\special{pa 1600 1460}%
\special{pa 1470 1580}%
\special{fp}%
\special{pa 1630 1260}%
\special{pa 1040 400}%
\special{fp}%
\special{pa 1050 550}%
\special{pa 1050 400}%
\special{fp}%
\special{pa 1050 400}%
\special{pa 1200 450}%
\special{fp}%
\put(10.0000,-2.2000){\makebox(0,0){1}}%
\put(1.9000,-14.1000){\makebox(0,0){2}}%
\put(18.0000,-14.0000){\makebox(0,0){3}}%
\end{picture}%
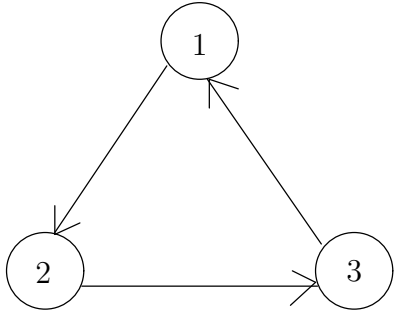
\\ Figure 4. $ \sigma^{(3)} $
}
\minip{7cm}{
\unitlength 0.1in
\begin{picture}( 20.0000, 16.0000)(  0.0000,-16.0000)
%
\special{pn 8}%
\special{ar 1000 200 200 200  0.0000000 6.2831853}%
%
\special{pn 8}%
\special{ar 200 1400 200 200  0.0000000 6.2831853}%
%
\special{pn 8}%
\special{ar 1800 1400 200 200  0.0000000 6.2831853}%
%
\special{pn 8}%
\special{pa 380 1300}%
\special{pa 950 400}%
\special{fp}%
\special{pa 800 470}%
\special{pa 940 400}%
\special{fp}%
\special{pa 940 400}%
\special{pa 940 580}%
\special{fp}%
\special{pa 1600 1320}%
\special{pa 400 1320}%
\special{fp}%
\special{pa 520 1230}%
\special{pa 390 1310}%
\special{fp}%
\special{pa 390 1310}%
\special{pa 530 1390}%
\special{fp}%
\special{pa 1170 320}%
\special{pa 1780 1200}%
\special{fp}%
\special{pa 1650 1150}%
\special{pa 1790 1200}%
\special{fp}%
\special{pa 1790 1200}%
\special{pa 1800 1070}%
\special{fp}%
\put(10.0000,-2.2000){\makebox(0,0){1}}%
\put(1.9000,-14.1000){\makebox(0,0){2}}%
\put(18.0000,-14.0000){\makebox(0,0){3}}%
\end{picture}%
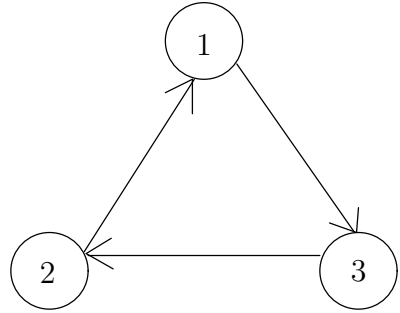
\\ Figure 5. $ \sigma^{(4)} $
}
\end{center}

\subsection{Realization of mixing Markov chain as a $ \mu $-random walk}

Our aim is to choose a mapping law $ \mu $ which satisfies a nice property. 

\begin{Def} \label{def: synchr}
A subset $ \Sigma_0 $ of $ \Sigma $ is called {\em synchronizing} 
if there exists a sequence $ s=(\sigma_p,\ldots,\sigma_1) $ of elements of $ \Sigma_0 $ 
such that the composition product $ \abra{s}:=\sigma_p \cdots \sigma_1 $ maps $ V $ 
onto a singleton. 
\end{Def}

First one of our main theorems is the following. 

\begin{Thm} \label{thm: main}
Suppose that $ Y=(Y_k)_{k \in \bZ} $ is mixing. 
Then one can choose a mapping law $ \mu $ for $ \Trans $ 
so that $ \mu $ has synchronizing support. 
\end{Thm}

Theorem \ref{thm: main} will be proved in Section \ref{sec: prf1}. 

Let us explain how our $ \mu $-random walk is related to road coloring. 
The support of $ \mu $, which we denote by $ \{ \sigma^{(1)},\ldots,\sigma^{(d)} \} $, 
induces the adjacency matrix $ A $ of a directed graph $ (V,A) $ 
which is of constant outdegree, i.e., 
from every site there are $ d $ roads laid. 
Then each element $ \sigma^{(1)},\ldots,\sigma^{(d)} $ may be regarded as a road color 
so that no two roads from the same site have the same color. 
For a $ \mu $-random walk $ (X,N) $, 
the process $ X $ moves in the directed graph $ (V,A) $ 
being driven by the randomly-chosen road colors indicated by $ N $ 
via equation \eqref{eq: Tsirel}. 
Thus we may call $ (X,N) $ 
{\em a random walk in a directed graph subject to a road coloring}. 
For example, the directed graphs induced by $ \mu^{(1)} $ and $ \mu^{(2)} $ 
which are defined in \eqref{eq: mu1} and \eqref{eq: mu2}, respectively, 
are illustrated as Figures 6 and 7, respectively. 
Since $ \sigma^{(1)} \sigma^{(2)} V = \{ 3 \} $, 
we see that the support of $ \mu^{(1)} $ is synchronizing, 
while we can easily see that the support of $ \mu^{(2)} $ is non-synchronizing. 

\begin{center}
\minip{7cm}{
\unitlength 0.1in
\begin{picture}( 20.0000, 16.0000)(  0.0000,-16.0000)
%
\special{pn 8}%
\special{ar 1000 200 200 200  0.0000000 6.2831853}%
%
\special{pn 8}%
\special{ar 200 1400 200 200  0.0000000 6.2831853}%
%
\special{pn 8}%
\special{ar 1800 1400 200 200  0.0000000 6.2831853}%
%
\special{pn 24}%
\special{pa 1170 320}%
\special{pa 1780 1200}%
\special{fp}%
\special{pa 1650 1150}%
\special{pa 1790 1200}%
\special{fp}%
\special{pa 1790 1200}%
\special{pa 1800 1070}%
\special{fp}%
\special{pa 390 1480}%
\special{pa 1610 1480}%
\special{fp}%
\special{pa 1480 1400}%
\special{pa 1600 1460}%
\special{fp}%
\special{pa 1600 1460}%
\special{pa 1470 1580}%
\special{fp}%
\special{pa 1630 1260}%
\special{pa 1040 400}%
\special{fp}%
\special{pa 1050 550}%
\special{pa 1050 400}%
\special{fp}%
\special{pa 1050 400}%
\special{pa 1200 450}%
\special{fp}%

\special{pn 8}%
\special{pa 830 330}%
\special{pa 250 1200}%
\special{fp}%
\special{pa 250 1060}%
\special{pa 250 1210}%
\special{fp}%
\special{pa 250 1210}%
\special{pa 380 1150}%
\special{fp}%
\special{pa 380 1300}%
\special{pa 950 400}%
\special{fp}%
\special{pa 800 470}%
\special{pa 940 400}%
\special{fp}%
\special{pa 940 400}%
\special{pa 940 580}%
\special{fp}%
\special{pa 1600 1320}%
\special{pa 400 1320}%
\special{fp}%
\special{pa 520 1230}%
\special{pa 390 1310}%
\special{fp}%
\special{pa 390 1310}%
\special{pa 530 1390}%
\special{fp}%
%
\special{pn 14}%
\special{ar 1440 1320 1316 1316  3.2228945 4.2118454}%
%
\special{pa 40 1100}%
\special{pa 120 1230}%
\special{fp}%
\special{pa 110 1210}%
\special{pa 200 1110}%
\special{fp}%
%
\special{ar 710 1310 1230 1230  5.1144264 6.2321348}%
%
\special{pa 1320 140}%
\special{pa 1190 180}%
\special{fp}%
\special{pa 1180 180}%
\special{pa 1300 350}%
\special{fp}%
%
\special{ar 980 200 1550 1550  1.0829599 2.0344439}%
%
\special{pa 1530 1580}%
\special{pa 1680 1570}%
\special{fp}%
\special{pa 1680 1570}%
\special{pa 1630 1710}%
\special{fp}%
\put(10.0000,-2.2000){\makebox(0,0){1}}%
\put(1.9000,-14.1000){\makebox(0,0){2}}%
\put(18.0000,-14.0000){\makebox(0,0){3}}%
\put(4.6000,-6.5000){\makebox(0,0){$ \sigma^{(2)} $}}%
\put(15.8000,-6.5000){\makebox(0,0){$ \sigma^{(1)} $}}%
\put(7.5000,-9.0000){\makebox(0,0){$ \sigma^{(2)} $}}%
\put(12.5000,-9.0000){\makebox(0,0){$ \sigma^{(1)} $}}%
\put(10.0000,-12.2000){\makebox(0,0){$ \sigma^{(2)} $}}%
\put(9.0000,-15.6000){\makebox(0,0){$ \sigma^{(1)} $}}%
\put(12.4000,-16.7000){\makebox(0,0){$ \sigma^{(3)} $}}%
\put(4.6000,-2.8000){\makebox(0,0){$ \sigma^{(3)} $}}%
\put(15.8000,-2.5000){\makebox(0,0){$ \sigma^{(3)} $}}%
\end{picture}%
\\[3mm] Figure 6. The graph induced by $ \mu^{(1)} $
}
\minip{7cm}{
\unitlength 0.1in
\begin{picture}( 20.0000, 16.0000)(  0.0000,-16.0000)
%
\special{pn 8}%
\special{ar 1000 200 200 200  0.0000000 6.2831853}%
%
\special{pn 8}%
\special{ar 200 1400 200 200  0.0000000 6.2831853}%
%
\special{pn 8}%
\special{ar 1800 1400 200 200  0.0000000 6.2831853}%
%
\special{pn 20}%
\special{pa 830 330}%
\special{pa 250 1200}%
\special{fp}%
\special{pa 250 1060}%
\special{pa 250 1210}%
\special{fp}%
\special{pa 250 1210}%
\special{pa 380 1150}%
\special{fp}%
\special{pa 390 1480}%
\special{pa 1610 1480}%
\special{fp}%
\special{pa 1480 1400}%
\special{pa 1600 1460}%
\special{fp}%
\special{pa 1600 1460}%
\special{pa 1470 1580}%
\special{fp}%
\special{pa 1630 1260}%
\special{pa 1040 400}%
\special{fp}%
\special{pa 1050 550}%
\special{pa 1050 400}%
\special{fp}%
\special{pa 1050 400}%
\special{pa 1200 450}%
\special{fp}%
\special{pn 8}%
\special{pa 380 1300}%
\special{pa 950 400}%
\special{fp}%
\special{pa 800 470}%
\special{pa 940 400}%
\special{fp}%
\special{pa 940 400}%
\special{pa 940 580}%
\special{fp}%
\special{pa 1600 1320}%
\special{pa 400 1320}%
\special{fp}%
\special{pa 520 1230}%
\special{pa 390 1310}%
\special{fp}%
\special{pa 390 1310}%
\special{pa 530 1390}%
\special{fp}%
\special{pa 1170 320}%
\special{pa 1780 1200}%
\special{fp}%
\special{pa 1650 1150}%
\special{pa 1790 1200}%
\special{fp}%
\special{pa 1790 1200}%
\special{pa 1800 1070}%
\special{fp}%
\put(10.0000,-2.2000){\makebox(0,0){1}}%
\put(1.9000,-14.1000){\makebox(0,0){2}}%
\put(18.0000,-14.0000){\makebox(0,0){3}}%
\put(4.6000,-6.5000){\makebox(0,0){$ \sigma^{(3)} $}}%
\put(15.8000,-6.5000){\makebox(0,0){$ \sigma^{(4)} $}}%
\put(7.5000,-9.0000){\makebox(0,0){$ \sigma^{(4)} $}}%
\put(12.5000,-9.0000){\makebox(0,0){$ \sigma^{(3)} $}}%
\put(10.0000,-12.2000){\makebox(0,0){$ \sigma^{(4)} $}}%
\put(10.0000,-16.0000){\makebox(0,0){$ \sigma^{(3)} $}}%
\end{picture}%
\\[2.3mm] Figure 7. The graph induced by $ \mu^{(2)} $
}
\end{center}

Let us come back to the general discussion. 
If $ (X,N) $ is a $ \mu $-random walk 
and if the support of $ \mu $ is synchronizing, 
then the process $ X $ may be represented as 
\begin{align}
X_k = F(N_k,N_{k-1},\ldots) 
, \quad k \in \bZ 
\label{eq: repre}
\end{align}
for some measurable function $ F:\Sigma^{-\bN} \to V $. 
In fact, define 
\begin{align}
T(N_j:j \le k) 
= \max \{ l \in \bZ, \ l<k : N_k N_{k-1} \cdots N_{l} V \ \text{is a singleton} \} , 
\label{}
\end{align}
where we follow the convention that $ \max \emptyset = - \infty $. 
Since the support of $ \mu $ is synchronizing, 
it holds that $ T(N_k,N_{k-1},\ldots) $ is finite a.s. for all $ k \in \bZ $, 
so that we may define 
\begin{align}
X_k = N_k N_{k-1} \cdots N_{T(N_j:j \le k)} x_0 
, \quad k \in \bZ 
\label{eq: stopped}
\end{align}
for a fixed element $ x_0 \in V $, 
but the resulting random walk does not depend on the choice of $ x_0 $. 
This is such a representation as \eqref{eq: repre}. 

Letting $ k=0 $ in the identity \eqref{eq: stopped}, we have 
\begin{align}
X_0 = N_0 N_{-1} \cdots N_{T(N_j:j \le 0)} x_0 . 
\label{eq: stopped 0}
\end{align}
This shows that the stationary law of the Markov chain 
may be simulated exactly from an IID sequence. 
This method was a central idea of 
{\em Propp--Wilson's coupling from the past} (\cite{MR1611693}; see also \cite{MR1908939}). 
Our Theorem \ref{thm: main} assures theoretically that 
for any mixing Markov chain there always exists an appropriate mapping law 
such that Propp--Wilson's algorithm terminates almost surely. 

For the study of $ \mu $-random walks in the case of non-synchronizing supports, 
see Yano \cite{Yrcp}. 
Equation \eqref{eq: Tsirel} is called {\em Tsirelson's equation in discrete time}; 
see Yor \cite{MR1147613}, Akahori--Uenishi--Yano \cite{MR2365485}, 
Yano--Takahashi \cite{YT}, 
Yano--Yor \cite{YY} and Hirayama--Yano \cite{HY,HYtsirel} for the details. 

The representation $ Y \dist X = F(N) $ of $ Y $ by an IID sequence $ N $ 
of the form \eqref{eq: repre} 
is called a {\em non-anticipating representation}. 
Rosenblatt (\cite{MR0114249} and \cite{MR0166839}) obtained a necessary and sufficient condition 
for a Markov chain with countable state space to have a non-anticipating representation 
$ Y \dist X = F(N) $ where $ N=(N_k)_{k \in \bZ} $ is an IID with uniform law on $ [0,1] $.

\subsection{Condition for approximate preservation of entropies}

Let $ Y $ and $ (X,N) $ as in Theorem \ref{thm: main}. 
Let us compare the amounts of information of $ Y $ and $ N $ in terms of their entropies. 
See the standard textbook \cite{MR524567} for basic theory of entropies. 
Let $ \lambda $ be the stationary law of $ Y $ and define 
\begin{align}
h(Y) = - \sum_{x,y \in V} \lambda(x) \trans_{x,y} \log \trans_{x,y} 
\label{}
\end{align}
and 
\begin{align}
h(N) = - \sum_{\sigma \in \Sigma} \mu(\sigma) \log \mu(\sigma) . 
\label{}
\end{align}
Since $ Y \dist X $ and $ X $ is a measurable function of $ N $ as in \eqref{eq: stopped}, 
we have 
\begin{align}
h(Y) \le h(N) . 
\label{eq: entropy ineq}
\end{align}
Note that Ornstein--Friedman's theorem (\cite{MR0257322} and \cite{MR0274718}) asserts that 
two mixing Markov chains which have common entropy are isomorphic. 
By this theorem, we see that, if equality holds in \eqref{eq: entropy ineq}, 
then $ Y $ is isomorphic to $ N $. 
We do not have any general criterion on $ Y $ for existence of a mapping law 
such that $ Y $ is isomorphic to $ N $. 
We will give an example for non-existence in Section \ref{sec: cex}. 

We are interested in when we can take mapping laws 
such that the $ h(N) $ approximates the $ h(Y) $. 
Following Rosenblatt \cite{MR0114249}, we introduce the following: 

\begin{Def}
A stationary Markov chain $ Y $ is called {\em p-uniform} if 
there exist a probability law $ \nu $ on $ V $ 
and a family $ \{ \tau_x:x \in V \} $ of permutations of $ V $ 
such that 
\begin{align}
\trans_{x,y} = \nu(\tau_x(y)) 
, \quad x,y \in V . 
\label{eq: nu}
\end{align}
\end{Def}

The second one of our main theorems is the following. 

\begin{Thm} \label{thm: entropy}
Let $ Y $ be a mixing Markov chain. Then the following assertions are equivalent: 
\begin{enumerate}
\item 
There exists a sequence $ \{ \mu^{(n)}:n=1,2,\ldots \} $ 
of mapping laws for $ \Trans $ with synchronizing support such that 
the $ N^{(n)} $ corresponding to $ \mu^{(n)} $ satisfies 
\begin{align}
h(N^{(n)}) \to h(Y) 
\quad \text{as $ n \to \infty $}. 
\label{eq: i-2}
\end{align}
\item 
$ Y $ is p-uniform. 
\end{enumerate}
In particular, if $ h(N)=h(Y) $ holds for $ N $ 
corresponding to some mapping law $ \mu $ for $ Q $ 
with synchronizing support, then $ Y $ is necessarily p-uniform. 
\end{Thm}

Theorem \ref{thm: entropy} will be proved in Section \ref{sec: prf2}.

This paper is organized as follows. 
In Section \ref{sec: rcp}, we introduce several notations 
to state the road coloring problem. 
Sections \ref{sec: prf1} and \ref{sec: prf2} 
are devoted to the proofs of Theorems \ref{thm: main} and \ref{thm: entropy}, respectively. 
In Section \ref{sec: cex}, we give an example for Theorem \ref{thm: entropy}.

\section{Road colorings of a directed graph} \label{sec: rcp}

Let $ A=[A(y,x)]_{y,x \in V} $ be a $ V \times V $-matrix 
whose entries are non-negative integers. 
The pair $ (V,A) $ may be called a {\em directed graph}, 
where, for $ x,y \in V $, the value $ A(y,x) $ is regarded 
as the number of directed edges from $ x $ to $ y $. 
The set $ V $ is called {\em the set of vertices} 
and the matrix $ A $ is called {\em the adjacency matrix}. 

The graph $ (V,A) $ is called {\em of constant outdegree} if 
there exists a constant $ d $ such that 
\begin{align}
\sum_{y \in V} A(y,x) = d 
\quad \text{for all $ x \in V $}. 
\label{}
\end{align}
In this case $ (V,A) $ is called {\em $ d $-out}. 
The graph $ (V,A) $ is called {\em strongly connected} 
if, for any $ x,y \in V $, there exists a positive integer $ n=n(x,y) $ 
such that $ A^n(y,x) \ge 1 $. 
The graph $ (V,A) $ is called {\em aperiodic} 
if the greatest common divisor among $ \{ n \ge 1 : A^n(x,x) \ge 1 \} $ is one 
for all $ x \in V $. 
Note that $ (V,A) $ is both strongly connected and aperiodic 
if and only if there exists a positive integer $ r $ such that 
$ A^r(y,x) \ge 1 $ for all $ x,y \in V $. 
We say that the graph $ (V,A) $ or the adjacency matrix $ A $ 
{\em satisfies the assumption} {\bf (A)} if 
$ (V,A) $ is of constant outdegree, strongly connected and aperiodic.

Recall that $ \Sigma $ is the set of all mappings from $ V $ to itself. 
For $ \sigma_1,\sigma_2 \in \Sigma $ and $ x \in V $, 
we write $ \sigma_2 \sigma_1 x $ simply for $ \sigma_2 ( \sigma_1 (x) ) $. 
The set $ \Sigma $ acts $ V $ in the following sense: 
\begin{align}
(\sigma_1 \sigma_2) x 
=& \sigma_1 (\sigma_2 x) 
, \quad \sigma_1,\sigma_2 \in \Sigma , \ x \in V . 
\label{eq: action}
\end{align}
The set $ V = \{ 1,\ldots,m \} $ may be identified with the set of standard basis 
$ \{ e_1,\ldots,e_m \} $ of $ \bR^m $. 
An element $ \sigma \in \Sigma $ may be identified 
with the 1-out adjacency matrix $ \sigma=[\sigma(y,x)]_{y,x \in V} $ given as 
\begin{align}
\sigma = \bmat{ \sigma e_1 & \cdots & \sigma e_m } . 
\label{}
\end{align}
Under these identifications, we see that, for all $ x,y \in V $, 
\begin{align}
\sigma(y,x) = 1 
\quad \text{if and only if} \quad 
y = \sigma x . 
\label{eq: sigma and trans mat}
\end{align}

Let $ (V,A) $ be a $ d $-out directed graph. 
A family $ \{ \sigma^{(1)},\ldots,\sigma^{(d)} \} $ of elements of $ \Sigma $ 
(possibly with repeated elements) 
is called a {\em road coloring} of $ (V,A) $ if 
\begin{align}
A = \sigma^{(1)} + \cdots + \sigma^{(d)} . 
\label{eq: A}
\end{align}
Each $ \sigma^{(i)} $ is called a {\em road color}. 
Note that there exists at least one road coloring of $ (V,A) $. 
Conversely, 
if we are given a family $ \{ \sigma^{(1)},\ldots,\sigma^{(d)} \} $ of elements of $ \Sigma $ 
(possibly with repeated elements), 
then it induces a unique $ d $-out directed graph $ (V,A) $ given as \eqref{eq: A}. 

Let $ \Sigma_0 $ be a subset of $ \Sigma $. 
A sequence $ s = (\sigma_p,\ldots,\sigma_2,\sigma_1) $ of elements of $ \Sigma_0 $ 
is called a {\em $ \Sigma_0 $-word}. 
For a $ \Sigma_0 $-word $ s = (\sigma_p,\ldots,\sigma_2,\sigma_1) $, 
we write $ \abra{s} $ for the product $ \sigma_p \cdots \sigma_2 \sigma_1 $. 
The following definition is a slight modification of Definition \ref{def: synchr}. 

\begin{Def}
A road coloring $ \Sigma_0 = \{ \sigma^{(1)},\ldots,\sigma^{(d)} \} $ 
is called {\em synchronizing} 
if $ \Sigma_0 $ as a subset of $ \Sigma $ is synchronizing. 
\end{Def}

By this definition, we see that 
a road coloring $ \Sigma_0 = \{ \sigma^{(1)},\ldots,\sigma^{(d)} \} $ is synchronizing 
if and only if 
$ \abra{s} V $ is a singleton for some $ \Sigma_0 $-word $ s $. 
If we express 
\begin{align}
s = (\sigma^{(i(p))},\ldots,\sigma^{(i(2))},\sigma^{(i(1))}) 
\label{}
\end{align}
with some numbers $ i(1),\ldots,i(p) \in \{ 1,\ldots,d \} $, 
the assertion ``$ \abra{s} V $ is a singleton" may be stated in other words as follows: 
Those who walk in the graph $ (V,A) $ 
according to the colors $ \sigma^{(i(1))},\ldots,\sigma^{(i(p))} $ in this order 
will lead to a common vertex, no matter where they started from. 

Now we state the {\em road coloring theorem}. 

\begin{Thm}[{Trahtman (\cite{MR2534238})}] \label{thm: Trahtman}
Suppose that the directed graph $ (V,A) $ satisfies the assumption {\bf (A)}. 
Then there exists a synchronizing road coloring of $ (V,A) $. 
\end{Thm}

This was first conjectured in the case of no multiple directed edges 
by Adler--Goodwyn--Weiss \cite{MR0437715} (see also \cite{MR0257315}) 
in the context of the isomorphism problem of symbolic dynamics with common topological entropy. 
For related studies before Trahtman \cite{MR2534238}, 
see \cite{MR1788119}, \cite{MR2312963} and \cite{MR2081291}.

\section{Construction of a mapping law on a synchronizing road coloring} \label{sec: prf1}

We need the following lemma. 

\begin{Lem} \label{lem: rational}
Let $ Y $ be a stationary Markov chain with one-step transition probability matrix $ \Trans $. 
Then there exists a mapping law $ \mu $ for $ \Trans $. 
\end{Lem}

\Proof{
First, we suppose that $ \trans_{x,y} $ is a rational number for all $ x,y \in V $. 
Then we may take an integer $ d $ sufficiently large so that 
$ A(y,x) := \trans_{x,y} d $ is an integer for all $ x,y \in V $. 
Then $ A:=[A(y,x)]_{x,y \in V} $ is the adjacency matrix 
of a $ d $-out directed graph $ (V,A) $; 
in fact, 
\begin{align}
\sum_{y \in V} A(y,x) = d \sum_{y \in V} \trans_{x,y} = d . 
\label{}
\end{align}
Let $ \{ \sigma^{(1)},\ldots,\sigma^{(d)} \} $ be a road coloring of $ (V,A) $ 
and define 
\begin{align}
\mu(\sigma) = \frac{1}{d} \sharp (\{ i=1,\ldots,d: \sigma^{(i)}=\sigma \}) 
\label{}
\end{align}
where $ \sharp(\cdot) $ denotes the number of elements of the set indicated. 
Thus, for any $ x,y \in V $, we see that 
\begin{align}
\sum_{\sigma \in \Sigma: y=\sigma x } \mu(\sigma) 
= \frac{1}{d} \sharp (\{ i=1,\ldots,d: \sigma^{(i)}(y,x)=1 \}) 
= \frac{1}{d} A(y,x) = \trans_{x,y} , 
\label{}
\end{align}
which shows that $ \mu $ is a mapping law for $ \Trans $. 

Second, we consider the general case. 
Let us take a sequence $ \{ \Trans^{(n)}:n=1,2,\ldots \} $ 
of one-step transition probability matrices 
such that 
$ \trans^{(n)}_{x,y} $ is a rational number for all $ n $ and $ x,y \in V $ 
and that $ \trans^{(n)}_{x,y} \to \trans_{x,y} $ as $ n \to \infty $ for all $ x,y \in V $. 
Then for any $ n $ there exists a mapping law $ \mu^{(n)} $ for $ \Trans^{(n)} $. 
Since $ \Sigma $ is a finite set, 
we can choose some subsequence $ \{ \mu^{(n(k))}:k=1,2,\ldots \} $ 
and some probability law $ \mu $ on $ \Sigma $ 
such that $ \mu^{(n(k))}(\sigma) \to \mu(\sigma) $ as $ k \to \infty $. 
This shows that $ \mu $ is a mapping law for $ \Trans $. 
The proof is now complete. 
}

Now we proceed to prove Theorem \ref{thm: main}. 

\Proof[Proof of Theorem \ref{thm: main}]{
Let $ Q=(q_{x,y})_{x,y \in V} $ be the one-step transition probability matrix 
for a mixing Markov chain $ Y $. 

First, we take an adjacency matrix $ A $ 
which is of constant outdegree such that 
\begin{align}
A(y,x) 
\begin{cases}
\ge 1 & \text{if $ \trans_{x,y}>0 $}, \\ 
= 0 & \text{if $ \trans_{x,y}=0 $}. 
\end{cases}
\label{eq: A(y,x)}
\end{align}
For this, we introduce a subset $ V \times V $ defined by 
\begin{align}
E = \{ (x,y) \in V \times V : \trans_{x,y} > 0 \} . 
\label{}
\end{align}
For each $ x \in V $, we define the outdegree of $ E $ at $ x $ by 
\begin{align}
d(x) = \sharp \cbra{ (x,y) \in E : y \in V } 
\label{}
\end{align}
and write $ d = \max_{x \in V} d(x) $ for the maximum outdegree of $ E $. 
For each $ x \in V $, we may choose a site $ \sigma(x) \in V $ so that 
$ (x,\sigma(x)) \in E $. Now we may set 
\begin{align}
A(y,x) = 
\begin{cases}
d-d(x)+1 & \text{if $ y=\sigma(x) $}, \\
1 & \text{if $ y \neq \sigma(x) $ and $ (x,y) \in E $}, \\
0 & \text{otherwise}. 
\end{cases}
\label{}
\end{align}
Then this $ (A(y,x))_{x,y \in V} $ is as desired. 

Since $ Y $ is a mixing Markov chain, there exists a positive integer $ r $ such that 
$ \trans^r_{x,y} > 0 $ for all $ x,y \in V $. 
Hence we have $ A^r(y,x) \ge 1 $ for all $ x,y \in V $; 
in fact, there exists a path $ x=x_0,x_1,\ldots,x_n=y $ such that 
$ \trans_{x_{k-1},x_k}>0 $ for $ k=1,2,\ldots,n $, 
which implies that 
$ A(x_k,x_{k-1}) \ge 1 $ for $ k=1,2,\ldots,n $. 
Thus we see that $ (V,A) $ satisfies the assumption {\bf (A)}. 
Thus we may apply Theorem \ref{thm: Trahtman} 
to obtain a synchronizing road coloring 
$ \{ \sigma^{(1)},\ldots,\sigma^{(d)} \} $ of $ (V,A) $. 
Define 
\begin{align}
\hat{\mu}(\sigma) = \frac{1}{d} \sharp (\{ i=1,\ldots,d: \sigma^{(i)}=\sigma \}) 
, \quad \sigma \in \Sigma 
\label{}
\end{align}
and define 
\begin{align}
\hat{\trans}_{x,y} = \sum_{\sigma \in \Sigma: y=\sigma x} \hat{\mu}(\sigma) 
, \quad x,y \in V . 
\label{}
\end{align}
Then $ \hat{\mu} $ is a mapping law for $ \hat{\Trans} $ 
and has synchronizing support. 
We also note that 
\begin{align}
\hat{\trans}_{x,y} = 0 
\quad \text{if $ (x,y) \notin E $}. 
\label{}
\end{align}

Let 
\begin{align}
\eps = \min \{ \trans_{x,y}:(x,y) \in E \} > 0 . 
\label{}
\end{align}
If $ \eps=1 $, then we have $ Q = \hat{Q} $, 
so that $ \hat{\mu} $ is as desired. 
Let us assume that $ \eps<1 $. 
Define 
\begin{align}
\Trans^{(\eps)} = \frac{1}{1-\eps} \rbra{ \Trans - \eps \hat{\Trans} } . 
\label{}
\end{align}
Then $ \Trans^{(\eps)}=(\trans^{(\eps)}_{x,y})_{x,y \in V} $ 
is a one-step transition probability matrix 
of a stationary Markov chain. 
In fact, we see that 
\begin{align}
(1-\eps) \trans^{(\eps)}_{x,y} = \trans_{x,y} - \eps \hat{\trans}_{x,y} 
\ge \trans_{x,y} - \eps 1_{\{ (x,y) \in E \}} \ge 0 
, \quad x,y \in V 
\label{}
\end{align}
and that 
\begin{align}
\sum_{y \in V} \trans^{(\eps)}_{x,y} 
= \frac{1}{1-\eps} \rbra{ \sum_{y \in V} \trans_{x,y} - \eps \sum_{y \in V} \hat{\trans}_{x,y} } 
= 1. 
\label{}
\end{align}
Now we apply Lemma \ref{lem: rational} to obtain a mapping law 
$ \mu^{(\eps)} $ for $ \Trans^{(\eps)} $. 
Define 
\begin{align}
\mu = (1-\eps) \mu^{(\eps)} + \eps \hat{\mu} . 
\label{}
\end{align}
Since $ \mu^{(\eps)} $ has synchronizing support, 
so does $ \mu $. 
For $ x,y \in V $, we have 
\begin{align}
\sum_{\sigma \in \Sigma : y=\sigma x} \mu(\sigma) 
=& (1-\eps) \sum_{\sigma \in \Sigma : y=\sigma x} \mu^{(\eps)}(\sigma) 
+ \eps \sum_{\sigma \in \Sigma : y=\sigma x} \hat{\mu}(\sigma) . 
\label{} \\
=& (1-\eps) \trans^{(\eps)}_{x,y} + \eps \hat{\trans}_{x,y} 
= \trans_{x,y} , 
\label{}
\end{align}
which shows that $ \mu $ is a mapping law for $ \Trans $. 
The proof is now complete. 
}

\section{Approximate preservation of entropies} \label{sec: prf2}

Let us prove Theorem \ref{thm: entropy}. 

\Proof[Proof of Theorem \ref{thm: entropy}]{
Let us prove that (i) implies (ii). 
Note that 
\begin{align}
h(Y) =& - \sum_{x,y \in V} \lambda(x) \trans_{x,y} \log \trans_{x,y} 
, \label{eq: h(Y)} \\
h(N^{(n)}) =& - \sum_{\sigma \in \Sigma} \mu^{(n)}(\sigma) \log \mu^{(n)}(\sigma) . 
\label{eq: h(N)}
\end{align}
Taking a subsequence if necessary, we may assume that 
there exists a probability law $ \mu $ on $ \Sigma $ such that 
$ \mu^{(n)}(\sigma) \to \mu(\sigma) $ for all $ \sigma \in \Sigma $. 
Note that $ \mu $ is a mapping law for $ \Trans $ 
but does not necessarily have synchronizing support. 
By the assumption \eqref{eq: i-2}, we see that 
\begin{align}
h(Y) 
= \lim_{n \to \infty } h(N^{(n)}) 
= - \sum_{\sigma \in \Sigma} \mu(\sigma) \log \mu(\sigma) . 
\label{}
\end{align}

For $ x,y \in V $, we set 
\begin{align}
\Sigma(y,x) = \cbra{ \sigma \in \Sigma: y = \sigma x } , 
\label{}
\end{align}
so that we have 
\begin{align}
\trans_{x,y} = \sum_{\sigma \in \Sigma(y,x)} \mu(\sigma) . 
\label{}
\end{align}
Hence we have 
\begin{align}
\mu(\sigma) \le \trans_{x,y} 
\quad \text{whenever} \ \sigma \in \Sigma(y,x). 
\label{eq: mu and p}
\end{align}
Since $ t \mapsto \log t $ is increasing, we have 
\begin{align}
- \sum_{\sigma \in \Sigma(y,x)} \mu(\sigma) \log \mu(\sigma) 
\ge - \sum_{\sigma \in \Sigma(y,x)} \mu(\sigma) \log \trans_{x,y} 
= - \trans_{x,y} \log \trans_{x,y} . 
\label{eq: entr ineq}
\end{align}
Since $ \bigcup_{y \in V} \Sigma(y,x) = \Sigma $, we have 
\begin{align}
h(Y) = 
- \sum_{y \in V} \sum_{\sigma \in \Sigma(y,x)} \mu(\sigma) \log \mu(\sigma) 
\ge q(x) 
\quad \text{for all $ x \in V $} , 
\label{eq: h(N) ineq}
\end{align}
where we set 
\begin{align}
q(x) = - \sum_{y \in V} \trans_{x,y} \log \trans_{x,y} 
, \quad x \in V . 
\label{}
\end{align}

We take $ \hat{x} \in V $ such that 
\begin{align}
q(\hat{x}) = \max_{x \in V} q(x) . 
\label{}
\end{align}
Using \eqref{eq: h(N) ineq} and \eqref{eq: h(Y)}, we have 
\begin{align}
q(\hat{x}) \le h(Y) = \sum_{x \in V} \lambda(x) q(x) \le q(\hat{x}) . 
\label{eq: qhatx}
\end{align}
Thus we see that the equalities hold in \eqref{eq: qhatx} 
and that $ q(x)=q(\hat{x}) $ for all $ x \in V $. 
For any $ x \in V $, 
we combine $ h(N)=q(x) $ together with \eqref{eq: entr ineq} 
and then obtain 
\begin{align}
- \sum_{\sigma \in \Sigma(y,x)} \mu(\sigma) \log \mu(\sigma) 
= - \trans_{x,y} \log \trans_{x,y} 
, \quad x,y \in V . 
\label{}
\end{align}
Combining this with \eqref{eq: mu and p}, we obtain 
\begin{align}
\mu(\sigma) = \trans_{x,y} 
\quad \text{whenever} \ \sigma \in \Sigma(y,x). 
\label{}
\end{align}

Let $ x_0 \in V $ be fixed and let $ x \in V $. 
Since $ \{ \Sigma(y,x) : y \in V \} $ is a partition of $ \Sigma $, 
we may choose a permutation $ \tau_x $ of $ V $ so that 
\begin{align}
\Sigma(\tau_x(y),x) \cap \Sigma(y,x_0) \neq \emptyset 
, \quad y \in V . 
\label{}
\end{align}
This shows that 
\begin{align}
\trans_{x,\tau_x(y)} = \trans_{x_0,y} 
, \quad x,y \in V , 
\label{}
\end{align}
which implies p-uniformity of $ Y $. The proof is now complete. 

Let us prove that (ii) implies (i). 
Let $ \{ x_1,\ldots,x_d \} $ be an enumeration of 
the support of the law $ \nu $ in \eqref{eq: nu}. 
For $ i=1,\ldots,d $, we define 
\begin{align}
\sigma^{(i)}(y,x) = 1_{\{ \tau_x(y) = x_i \}} . 
\label{}
\end{align}
For each $ x \in V $, there exists a unique $ y \in V $ such that 
$ \sigma^{(i)}(y,x) = 1 $, so that we have $ \sigma^{(i)} \in \Sigma $. 
By \eqref{eq: nu}, we obtain 
\begin{align}
\trans_{x,y} = \sum_{i=1}^d \sigma^{(i)}(y,x) \nu(x_i) 
, \quad x,y \in V . 
\label{}
\end{align}

Let $ A $ be as in \eqref{eq: A(y,x)} 
and let $ \Sigma_1 $ be a synchronizing subset 
corresponding to some synchronizing road coloring of $ (V,A) $. 
For sufficiently large integer $ n $, 
we define a probability law $ \mu^{(n)} $ on $ \Sigma $ by 
\begin{align}
\mu^{(n)}(\sigma) 
= \sum_{i:\sigma^{(i)} = \sigma} \cbra{ \nu(x_i) - \frac{1}{nd} } 
+ \frac{1}{n |\Sigma_1|} 1_{\{ \sigma \in \Sigma_1 \}} . 
\label{}
\end{align}
Then it is obvious that $ \mu^{(n)} $ is a mapping law for $ \Trans $ 
and has synchronizing support. 

Let us verify the condition \eqref{eq: i-2}. 
On one hand, we have 
\begin{align}
h(N^{(n)}) \tend{}{n \to \infty } 
- \sum_{i=1}^d \nu(x_i) \log \nu(x_i) . 
\label{}
\end{align}
On the other hand, we have 
\begin{align}
h(Y) 
=& - \sum_{x,y \in V} \lambda(x) \trans_{x,y} \log \trans_{x,y} 
\label{} \\
=& - \sum_{x,y \in V} \lambda(x) \sum_{i=1}^d \sigma^{(i)}(y,x) \nu(x_i) \log \nu(x_i) 
\label{} \\
=& - \sum_{i=1}^d \cbra{ \sum_{x,y \in V} \lambda(x) \sigma^{(i)}(y,x) } \nu(x_i) \log \nu(x_i) 
\label{} \\
=& - \sum_{i=1}^d \nu(x_i) \log \nu(x_i) . 
\label{}
\end{align}
This shows \eqref{eq: i-2}. 
The proof is now complete. 
}

\section{An example} \label{sec: cex}

Let $ V=\{ 1,2 \} $. Then $ \Sigma = \{ (12),(21),(11),(22) \} $ where 
\begin{align}
(ij) 
= \bmat{ 1 \mapsto i \\ 2 \mapsto j } 
, \quad i,j=1,2. 
\label{}
\end{align}
Let $ 0<p<1 $ and consider a Markov chain $ Y $ with one-step transition probability given by 
\begin{align}
\bmat{ \trans_{1,1} & \trans_{1,2} \\ \trans_{2,1} & \trans_{2,2} } 
= 
\bmat{ p & 1-p \\ 1-p & p } . 
\label{}
\end{align}
Then it is obvious that $ Y $ is a mixing Markov chain. 
Since 
\begin{align}
\bmat{ \trans_{1,1} \\ \trans_{2,1} } 
= \bmat{ \trans_{2,2} \\ \trans_{1,2} } 
= \bmat{ p \\ 1-p } , 
\label{}
\end{align}
we see that $ Y $ is p-uniform. 

It is obvious that the stationary law is given as 
\begin{align}
\lambda(1) = \lambda(2) = \frac{1}{2}. 
\label{}
\end{align}
We now see that 
\begin{align}
h(Y) = \varphi(p) + \varphi(1-p) 
\label{}
\end{align}
where $ \varphi(t) = - t \log t $. 

If $ \mu $ is a mapping law for $ \Trans $, then we have 
\begin{align}
\mu(12) + \mu(11) = p 
, \quad 
\mu(21) + \mu(11) = 1-p . 
\label{}
\end{align}
From this, we see that there exists some $ \eps $ with $ 0 \le \eps \le \min \{ p,1-p \} $ 
such that 
\begin{align}
\eps = \mu(11) = \mu(22) 
, \quad 
\mu(12) = p - \eps 
, \quad 
\mu(21) = 1 - p - \eps . 
\label{eq: cex mu}
\end{align}
Conversely, for any $ \eps $ with $ 0 \le \eps \le \min \{ p,1-p \} $, 
we may define $ \mu=\mu^{(\eps)} $ by equation \eqref{eq: cex mu} 
so that $ \mu^{(\eps)} $ is a mapping law for $ \Trans $. 

If $ \mu^{(\eps)} $ has synchronizing support, $ \eps $ should be positive. 
Let $ \{ X^{(\eps)},N^{(\eps)} \} $ be the $ \mu^{(\eps)} $-random walk. 
We then see that 
\begin{align}
h(N^{(\eps)}) = 2 \varphi(\eps) + \varphi(p-\eps) + \varphi(1-p-\eps) . 
\label{}
\end{align}

If $ p=1/2 $, we see that $ h(Y)=h(N^{(1/2)}) $. 

Suppose that $ p \neq 1/2 $. 
Then, by an easy computation, we may see that 
\begin{align}
h(Y) < h(N^{(\eps)}) 
\label{}
\end{align}
for all $ \eps $ with $ 0 < \eps \le \min \{ p,1-p \} $. 
However, it holds that $ h(N^{(\eps)}) \to h(Y) $ as $ \eps \to 0+ $.

\def\cprime{$'$} \def\cprime{$'$}


\begin{thebibliography}{10}

\bibitem{MR0437715}
R.~L. Adler, L.~W. Goodwyn, and B.~Weiss.
\newblock Equivalence of topological {M}arkov shifts.
\newblock {\em Israel J. Math.}, 27(1):48--63, 1977.

\bibitem{MR0257315}
R.~L. Adler and B.~Weiss.
\newblock {\em Similarity of automorphisms of the torus}.
\newblock Memoirs of the American Mathematical Society, No. 98. American
  Mathematical Society, Providence, R.I., 1970.

\bibitem{MR2365485}
J.~Akahori, C.~Uenishi, and K.~Yano.
\newblock Stochastic equations on compact groups in discrete negative time.
\newblock {\em Probab. Theory Related Fields}, 140(3-4):569--593, 2008.

\bibitem{MR524567}
P.~Billingsley.
\newblock {\em Ergodic theory and information}.
\newblock Robert E. Krieger Publishing Co., Huntington, N.Y., 1978.
\newblock Reprint of the 1965 original.

\bibitem{MR2081291}
G.~Budzban.
\newblock Semigroups and the generalized road coloring problem.
\newblock {\em Semigroup Forum}, 69(2):201--208, 2004.

\bibitem{MR2312963}
G.~Budzban and P.~Feinsilver.
\newblock Completely simple semigroups, {L}ie algebras, and the road coloring
  problem.
\newblock {\em Semigroup Forum}, 74(2):206--226, 2007.

\bibitem{MR1788119}
G.~Budzban and A.~Mukherjea.
\newblock A semigroup approach to the road coloring problem.
\newblock In {\em Probability on algebraic structures ({G}ainesville, {FL},
  1999)}, volume 261 of {\em Contemp. Math.}, pages 195--207. Amer. Math. Soc.,
  Providence, RI, 2000.

\bibitem{MR0274718}
N.~A. Friedman and D.~S. Ornstein.
\newblock On isomorphism of weak {B}ernoulli transformations.
\newblock {\em Advances in Math.}, 5:365--394 (1970), 1970.

\bibitem{MR1908939}
O.~H{\"a}ggstr{\"o}m.
\newblock {\em Finite {M}arkov chains and algorithmic applications}, volume~52
  of {\em London Mathematical Society Student Texts}.
\newblock Cambridge University Press, Cambridge, 2002.

\bibitem{HY}
T.~Hirayama and K.~Yano.
\newblock Extremal solutions for stochastic equations indexed by negative
  integers and taking values in compact groups.
\newblock {\em Stochastic Process. Appl.}, 120(8):1404--1423, 2010.

\bibitem{HYtsirel}
T.~Hirayama and K.~Yano.
\newblock Strong solutions of {T}sirelson's equation in discrete time taking
  values in compact spaces with semigroup action.
\newblock {\em Preprint, ar{X}iv:1005.0038, submitted}, 2010.

\bibitem{MR0257322}
D.~Ornstein.
\newblock Bernoulli shifts with the same entropy are isomorphic.
\newblock {\em Advances in Math.}, 4:337--352 (1970), 1970.

\bibitem{MR1611693}
J.~G. Propp and D.~B. Wilson.
\newblock Exact sampling with coupled {M}arkov chains and applications to
  statistical mechanics.
\newblock In {\em Proceedings of the {S}eventh {I}nternational {C}onference on
  {R}andom {S}tructures and {A}lgorithms ({A}tlanta, {GA}, 1995)}, volume~9,
  pages 223--252, 1996.

\bibitem{MR0114249}
M.~Rosenblatt.
\newblock Stationary processes as shifts of functions of independent random
  variables.
\newblock {\em J. Math. Mech.}, 8:665--681, 1959.

\bibitem{MR0166839}
M.~Rosenblatt.
\newblock Stationary {M}arkov chains and independent random variables.
\newblock {\em J. Math. Mech.}, 9:945--949, 1960.
\newblock Addendum: {\it ibid.} vol. 11, page 317, 1962.

\bibitem{MR2534238}
A.~N. Trahtman.
\newblock The road coloring problem.
\newblock {\em Israel J. Math.}, 172:51--60, 2009.

\bibitem{Yrcp}
K.~Yano.
\newblock Random walk in a finite directed graph subject to a road coloring.
\newblock {\em Preprint, ar{X}iv:1005.0079}, 2010.

\bibitem{YT}
K.~Yano and Y.~Takahashi.
\newblock Time evolution with and without remote past.
\newblock {\em S\=urikaisekikenky\=usho K\=oky\=uroku}, 1552:164--171, 2007.
\newblock Recent Developments in Dynamical Systems (Kyoto, 2006).

\bibitem{YY}
K.~Yano and M.~Yor.
\newblock Around {T}sirelson's equation, or: The evolution process may not
  explain everything.
\newblock {\em Preprint, ar{X}iv:0906.3442, submitted}, 2009.

\bibitem{MR1147613}
M.~Yor.
\newblock Tsirel\cprime son's equation in discrete time.
\newblock {\em Probab. Theory Related Fields}, 91(2):135--152, 1992.

\end{thebibliography}
\end{document}